\documentclass{article}
\usepackage[utf8]{inputenc}
\usepackage{amsmath, amssymb, mathtools}  
\usepackage{hyperref}
\usepackage{graphicx, subfigure}
\usepackage{caption}
\usepackage{todonotes}
\usepackage{bbm}
\usepackage{caption}
\usepackage{hyperref}

\def\hmath$#1${\texorpdfstring{{\rmfamily\textit{#1}}}{#1}}

\newcommand{\matr}[1]{\mathbf{#1}}
\newcommand{\Amat}{\matr{A}}
\newcommand{\vect}[1]{\boldsymbol{#1}}

\newcommand{\R}{\mathbb{R}}
\newcommand{\C}{\mathbb{C}}
\newcommand{\infnorm}[1]{\lVert#1\rVert_\infty}

\title{The infinity norm bounds and characteristic polynomial for high order RK matrices}
    \author{Gayatri \v{C}aklovi\'{c}}
\date{January 2022}

\begin{document}

\maketitle
\begin{abstract}
This paper shows that $t_m \leq \infnorm{\Amat} \leq \sqrt{t_m}$ holds, when $\Amat \in \R^{m \times m}$ is a Runge-Kutta matrix which nodes originating from the Gaussian quadrature that integrates polynomials of degree $2m-2$ exactly. 
It can be shown that this is also true for the Gauss-Lobatto quadrature.
Additionally, the characteristic polynomial of $\Amat$, when the matrix is nonsingular, is $p_A(\lambda) = m!t^m + (m-1)!a_{m-1}t^{m-1} + \dots + a_0$, where the coefficients $a_i$ are the coefficients of the polynomial of nodes $\omega(t) = (t - t_1) \dots (t - t_m) = t^m + a_{m-1}t^{m-1} + \dots + a_0$.
\end{abstract}

\section{The Runge-Kutta method and Gaussian quadrature nodes}
The Runge-Kutta methods are a family of iterative methods for approximating the solutions of ordinary differential equations. 
A method can be represented in a form of a famous Butcher tableau: 

\[\renewcommand\arraystretch{1.2}
\begin{array}
{c|c}
\vect{c} & \Amat \\
\hline
& \vect{b}^T
\end{array}, \]
where $\vect{c}, \vect{b} \in \R^m$ represent the nodes and the weights and the matrix $\Amat \in \R^{m \times m}$ is called the Runge-Kutta matrix of order $m$. 
More about the origin of the methods can be found in \cite{ODERK}.

\subsection{Gaussian quadrature nodes}
On function spaces, there is a natural scalar product. 
For integrable functions $f, g:\R \rightarrow \R$ a scalar product can be defined as
\begin{equation} \label{eq:scalarprod}
\langle f, g \rangle = \int_0^1f \,g .
\end{equation}
This allows us to define orthogonality on function spaces.
The family of orthogonal polynomials can be used to construct a high order Gaussian-quadrature rules.
Let
\begin{equation}\label{eq:quadrature}
    \int_0^1 p(s)ds = \sum_{i=1}^{m} w_i p(t_i)
\end{equation}
define a quadrature rule, where some points $t_i$ may be fixed and while other are roots of an orthogonal polynomial. Then, the right sum computes the integral correctly
for polynomials of degree $2m-2$ or less. 
Depending on the choice of these points, we get different quadrature rules such as Legendre, Radau, Lobatto, etc.
The integration weights $w_i$ are positive and the integration points $t_i$ are inside the integration interval \cite{QUAD}. 
In our case, without loss of generality, we will assume that $t_i \in [0, 1]$.

Let $l_j$ be a Lagrange interpolation polynomial defined in points $t_i$:
\[ l_j(t) = \prod_{i \neq j} \frac{t - t_i}{t_j - t_i}, \quad j = 1, \dots, m. \]
If we define the Runge-Kutta matrix as 
\begin{equation}\label{eq:highRK}
    \Amat_{ij} = \int_0^{t_i}l_j(s)ds,
\end{equation}
we end up with a high order Runge-Kutta method where the vector $\vect{c} = (t_1, \dots, t_m)$ and $\vect{b} = (w_1, \dots, w_m)$.

\section{Bounds for \texorpdfstring{$\|\Amat\|_\infty$}{}}
Let The Runge-Kutta matrix be defined as in \eqref{eq:highRK}. In order to compute lower and upper bounds for $\|\Amat\|_\infty$, we firstly have to understand what does a matrix-vector product $\Amat \vect{x}$ mean. The product behaves as 
\begin{equation*}
    \Amat \vect{x} = 
    \begin{bmatrix}
    \sum_{j = 1}^m \vect{x}_j \int_0^{t_1} l_j(s)ds \\
    \sum_{j = 1}^m \vect{x}_j \int_0^{t_2} l_j(s)ds \\
    \vdots \\
    \sum_{j = 1}^m \vect{x}_j \int_0^{t_m} l_j(s)ds \\
    \end{bmatrix}
    =
    \begin{bmatrix}
    \int_0^{t_1} \sum_{j = 1}^m \vect{x}_j l_j(s)ds \\
    \int_0^{t_2} \sum_{j = 1}^m \vect{x}_j l_j(s)ds \\
    \vdots \\
    \int_0^{t_m} \sum_{j = 1}^m \vect{x}_j l_j(s)ds \\
    \end{bmatrix}.
\end{equation*}
This motivates us to construct a one-to-one map from a vector $\vect{x}\in \R^m$ to a polynomial $p\in P_{m-1}$, where $p := \sum_{j = 1}^m \vect{x}_j l_j$.
This is uniquely defined because the polynomial $p$ passes through $m$ points $(t_i, \vect{x}_i)$ and is of degree $m-1$.
Then, the sum $\int_0^{t_i} \sum_{j = 1}^m \vect{x}_j l_j(s)ds$ is nothing more than $\int_0^{t_i}p(s)ds$ and we can conclude that multiplying a vector with a matrix $\Amat$ is as computing integrals of a corresponding polynomial.

Now, we are interested in computing the norm $\infnorm{\Amat}$. By definition, the norm is
\begin{equation}\label{eq:norm matrix}
    \infnorm{\Amat} = \max_{\infnorm{\vect{x}} = 1} \infnorm{\Amat \vect{x}},
\end{equation}
but the infinity norm of a matrix can also be computed as
\begin{equation*}
    \infnorm{\Amat} = \max_{1 \leq i \leq m} \sum_{j = 1}^m |\Amat_{ij}| = \max_{1 \leq i \leq m} \sum_{j = 1}^m \bigg|\int_0^{t_i} l_j(s)ds\bigg|,
\end{equation*}
which is the maximum row-sum of absolute values of its entries. From this, we can see that
\begin{align*}
    t_i = \bigg|\int_0^{t_i} 1ds \bigg| & =\bigg|\int_0^{t_i} \sum_{j = 1}^m l_j(s)ds \bigg| \\
    & = \bigg|\sum_{j = 1}^m \int_0^{t_i} l_j(s)ds \bigg| \\
    & \leq \sum_{j = 1}^m \bigg| \int_0^{t_i} l_j(s)ds \bigg| \\
    & \leq \sum_{j = 1}^m |\Amat_{ij}| \leq \infnorm{\Amat}.
\end{align*}
This yields a lower bound: $t_m \leq \infnorm{\Amat}$.

Now we will redefine the definition of the matrix norm $\eqref{eq:norm matrix}$ in a way that we do not compute the maximum over all vectors
$\vect{x}$ for which $\infnorm{x} = 1$, but over uniquely defined polynomials $p$ since $\R^m \ni \vect{x} \mapsto p \in P_{m-1}$ is a one-to-one map, as discussed above.
Proving $\infnorm{\Amat} \leq 1$ is then equivalent to proving
\begin{equation}\label{eq:new norm}
    \max_{1 \leq i \leq M} \bigg| \int_0^{t_i} p(s)ds \bigg| \leq 1
\end{equation}
holds, for $p \in P_{m-1}$ and $\max_{1 \leq i \leq M}|p(t_i)| = 1$.

\subsection{Exact quadrature for degree \texorpdfstring{$2m-2$}{} or higher}
First, using the integral Cauchy–Bunyakovsky-Schwarz inequality on the scalar product defined as in \eqref{eq:scalarprod}, we can get an upper bound on
$\big|\int_0^{t_i}p(s)ds\big|$ as
\begin{align}
    \bigg|\int_0^{t_i}p(s)ds\bigg|^2 = \bigg|\int_0^{t_i}p(s) \cdot 1 ds\bigg|^2 &\leq \bigg|\int_0^{t_i}p^2(s)ds\bigg| \bigg|\int_0^{t_i}1^2ds\bigg| \\
    & \leq t_i \int_0^{t_i}p^2(s)ds \\
    & \leq t_i \int_0^{1}p^2(s)ds.\label{eq:0t}
\end{align}
Since $p^2 \in P_{2m-2}$, using the quadrature rule \eqref{eq:quadrature} which correctly
computes integrals of polynomials of degree $2m-2$, we get
\begin{align*}
    \int_0^{1}p^2(s)ds = \sum_{i=1}^{M} w_i p^2(t_i) \leq \underbrace{\bigg(\sum_{i=1}^{M} w_i\bigg)}_{=1} \max_{1 \leq i \leq M} p^2(t_i) = \max_{1 \leq i \leq M} p^2(t_i).
\end{align*}
Because we are computing \eqref{eq:new norm} over polynomials which satisfy
$\max_{1 \leq i \leq M}|p(t_i)| = 1$, we know that $ \max_{1 \leq i \leq M} p^2(t_i) = 1$, thus 
\begin{equation}\label{eq:01}
    \int_0^{1}p^2(s)ds \leq \max_{1 \leq i \leq M} p^2(t_i) = 1.
\end{equation}
Now \eqref{eq:0t} + \eqref{eq:01} yields
\begin{equation*}
    \bigg|\int_0^{t_i}p(s)ds\bigg|^2 \leq t_i.
\end{equation*}
Thus,
\begin{equation*}
    t_i \leq \bigg|\int_0^{t_i}p(s)ds\bigg| \leq \sqrt{t_i}.
\end{equation*}

This proves that $t_m \leq \infnorm{\Amat} \leq \sqrt{t_m} \leq 1$ holds for nodes originating from the Gaussian quadrature that integrates polynomials of degree $2m-2$ exactly.

\subsection{Gauss-Lobatto nodes}
The Gauss-Lobatto quadrature integrates polynomials of accuracy $2m-3$, therefore the same
argument why is $\infnorm{\Amat} \leq 1$ does not hold. However, we will show that this is
still true for $m \geq 2$. The Gauss-Lobatto quadrature on $[0, 1]$ is 
defined as
\begin{equation}\label{eq:lobatto}
    \int_0^1p(s)ds = \sum_{i = 1}^m w_i p(t_i) + R_m,
\end{equation}
where $R_m = -c p^{(2m-2)}(\xi), \, c \geq 0, \, \xi \in [0, 1]$, see \cite{QUAD} for details. 
If we manage to bind $\int_0^1 p^2(s)ds \leq 1$, then using the same train of thought and arguments as in $\eqref{eq:0t}$, we are done.

Let $p^2(s) = as^{2m-2} + \dots$, where $a > 0$. We can write
\begin{equation*}
    \int_0^1 p^2(s)ds = \int_0^1 \big(p^2(s) - as^{2m-2}\big)ds + \int_0^1 as^{2m-2}ds.
\end{equation*}
The polynomial $p^2(s) - as^{2m-2}$ is of degree $2m-3$ and can be computed with
quadrature \eqref{eq:lobatto} as
\begin{equation*}
    \int_0^1 \big(p^2(s) - as^{2m-2}\big)ds = \sum_{i = 1}^m w_i \big(p^2(t_i) - at_i^{2m-2}\big).
\end{equation*}
This yields
\begin{align}
    \int_0^1 p^2(s)ds &= \sum_{i = 1}^m w_i \big(p^2(t_i) - at_i^{2m-2}\big) + \int_0^1 as^{2m-2}ds \\
    & \leq 1 - a\sum_{i = 1}^m w_i t_i^{2m-2} + \int_0^1 as^{2m-2}ds,\label{eq:lobatto error}
\end{align}
where the inequality holds because $\sum_{i = 1}^m w_i p^2(t_i) \leq 1$ for 
$\max_{1 \leq i \leq M}|p(t_i)| \leq 1$.
Now using the quadrature \eqref{eq:lobatto} on a function $t^{2m-2}$ yields
\begin{equation*}
    \int_0^1 s^{2m-2} = \sum_{i = 1}^m w_i t_i^{2m-2} - c(2m-2)!
\end{equation*}
and combining this with \eqref{eq:lobatto error} gives
\begin{equation*}
    \int_0^1 p^2(s)ds \leq 1 - ac(2m-2)! \leq 1,
\end{equation*}
because $ac \geq 0$.
This proves that $t_m \leq \infnorm{\Amat} \leq \sqrt{t_m} \leq 1$ is also true for the Gauss-Lobatto quadrature.

\section{The characteristic polynomial of \texorpdfstring{$\Amat$}{}}
Let 
$$\omega(t) = (t - t_1) \dots (t - t_m) = t^m + a_{m-1}t^{m-1} + \dots + a_0$$
be the polynomial of nodes. 
If $t_i$ are the roots of an orthogonal polynomial defining the Gaussian quadrature, then $\omega$ is exactly that scaled polynomial. We will show that the characteristic polynomial of $\Amat$ is
\begin{equation}\label{eq:pa}
    p_A(\lambda) = m!t^m + (m-1)!a_{m-1}t^{m-1} + \dots + a_0.
\end{equation}
The eigenvalue of matrix $\Amat$ is a scalar $\lambda \in \C$ that satisfies 
\[ \Amat \vect{x} = \lambda \vect{x},\]
for some $\vect{x} \neq 0$. 
In other words, using a one-on-one map to the polynomial space $P_{m-1}$, we can reformulate the eigenproblem as finding a polynomial $p \neq 0$ such that
\[ \int_0^{t_i} p(s)ds = \lambda p(t_i), \quad i = 1, \dots, m.\]
Substituting $p$ with $g'$, where $g\in P_m$ yields
\begin{equation}\label{eq:eig1}
    g(t_i) - g(0) = \lambda g'(t_i), \quad i = 1, \dots, m.
\end{equation}
From here we see that the eigenvector represented as a polynomial $g$ is not uniquely
defined, because if $g$ is a solution of \eqref{eq:eig1}, then $d_mg(t) + d_0$ is also a solution. Therefore, without loss of generality, we can assume that $g$ has the form
\[ g(t) = t^m + d_{m-1}t^{m-1} + \dots + d_1t. \]
Setting $d_m = 1$ is equivalent to imposing that the eigenvector is normalized which yields a uniqueness of the eigenvector. 
However, here we impose $\vect{x}$ has a fixed infinity norm, one that may deffer from 1. 
The part where we demand that $d_0 = 0$ comes from substituting $p$ with $g'$.
As a consequence we have $g(0) = 0$. 
With these assumptions, $g$ is uniquely defined since we have $m$ nonlinear equations and $m$ unknowns: $(d_{m-1}, \dots , d_1, \lambda)$ that are described with
\begin{equation}\label{eq:eig2}
    g(t_i) = \lambda g'(t_i), \quad i = 1, \dots, m.
\end{equation}
Let $G(t) := g(t) - \lambda g'(t)$, where g is the solution of \eqref{eq:eig2}. 
The difference $G - \omega$ is a polynomial of degree $m-1$ and both $\omega$ and $G$ are monic polynomials of degree $m$.
Because $G - \omega$ is zero in $m$ different points, we conclude that $G = \omega$.
Equations \eqref{eq:eig2} can now be rewritten in a continuous form as
\[ g(t) - \lambda g'(t) = \omega(t) .\]
Since the polynomials coefficients on the left side should be equal to the coefficients on the right side of the upper equation, we get
\begin{align*}
    -\lambda d_1 &= a_0 \\
    d_1 - 2 \lambda d_2 &= a_1 \\
    d_2 - 3 \lambda d_3 & = a_2 \\
    & \vdots \\
    d_{m-2} - (m-1)\lambda d_{m-1} &= a_{m-2}\\
    d_{m-1} - m\lambda &= a_{m-1}.
\end{align*}
Substituting the equations from the bottom to the second from the top, we get
\begin{equation}\label{eq:eig3}
    d_1 = m! \lambda^{m-1} + (m-1)!a_{m-1}\lambda^{m-2} + \dots + 2\lambda a_2,
\end{equation}
whereas from the first one we get $d_1 =- a_0 / \lambda$.
Substituting this in \eqref{eq:eig3} we end up with a polynomial $p_A$ as in \eqref{eq:pa} which zeros are the eigenvalues of our matrix $\Amat$. In other words, $p_A$ is exactly the characteristic polynomial of $\Amat$.

It remains to comment why $\lambda = 0$ is not a viable option for an eigenvalue. Matrix $\Amat$ is a mapping in a fashion
\[ (t_1^k, \dots, t_m^k) \rightarrow \frac{1}{k + 1}(t_1^{k + 1}, \dots, t_m^{k+1}), \quad k = 0, \dots, m-1, \]
if the quadrature rule defining $\Amat$ integrates polynomials of degree $m-1$ correctly.
The vectors $(t_1^k, \dots, t_m^k)$ form columns of the Vandermonde matrix $\vect{V}$ which is known to be nonsingular if the points $t_1, \dots, t_m$ are distinct. Because of this, we have $\Amat\vect{V} = \vect{D}_m\vect{V}$, where $\vect{D}_m = \operatorname{diag}(1, 1/2, \dots, 1/m)$. From here we see that $\Amat$ is nonsingular as a product of nonsingular matrices.

\section{Conclusion}
We proved that $t_m \leq \infnorm{\Amat} \leq \sqrt{t_m} \leq 1$ holds for nodes originating from the Gaussian quadrature that integrates polynomials of degree $2m-2$ exactly. 
Is is also true for the Gauss-Lobatto quadrature.
Additionally, the characteristic polynomial of $\Amat$ can be found, when the matrix is nonsingular.
The matrix is nonsingular if it originates from a quadrature that integrates polynomials in an exact way up to degree $m-1$.

\bibliography{lit}

\end{document}